\documentclass[draftcls,12pt,onecolumn]{IEEEtran}

\usepackage{times}
\usepackage{amsmath}
\usepackage{amsfonts}
\usepackage{amssymb}
\usepackage[pdftex]{graphicx}
\DeclareGraphicsExtensions{.pdf,.jpeg,.png}
\usepackage{cite}
\usepackage{url}
\usepackage{pdfpages}
\usepackage{epstopdf}
\usepackage{subcaption}

\title{\bf Permanent Magnet Synchronous Drives Observability Analysis for Motion-Sensorless Control}
\author{Mohamad Koteich, \textit{Student Member, IEEE}, \\ Abdelmalek Maloum, Gilles Duc and Guillaume Sandou
	\thanks{Mohamad Koteich is with Renault S.A.S. Technocentre, 78288 Guyancourt, France, and also with L2S - CentraleSup\'{e}lec - CNRS - Paris-Sud University, 91192 Gif-sur-Yvette, France (e-mail: mohamad.koteich@renault.com).}
	\thanks{Abdelmalek Maloum is with Renault S.A.S. Technocentre, 78288 Guyancourt, France (e-mail: abdelmalek.maloum@renault.com).}
	\thanks{Gilles Duc and Guillaume Sandou are with L2S - CentraleSup\'{e}lec - CNRS - Paris-Sud University, 91192 Gif-sur-Yvette, France  (e-mail: gilles.duc@centralesupelec.fr; guillaume.sandou@centralesupelec.fr).}}

\begin{document}
\maketitle
\thispagestyle{empty}
\pagestyle{empty}
\pagenumbering{arabic}
\section*{Keywords}
\begin{center}
Sensorless control, Synchronous motor, Induction Motor, AC machine.
\end{center}

\begin{abstract}
Motion-sensorless control techniques of electrical drives are attracting more attention in different industries. The local observability of sensorless permanent magnet synchronous drives is studied in this paper. A special interest is given to the standstill operation condition, where sensorless drives suffer of poor performance. Both, salient and non-salient machines are considered. The results are illustrated using numerical simulations.
\end{abstract}

\section{Introduction}

Permanent Magnet (PM) Synchronous Machine (SM) has been widely used in many potential industrial applications \cite{emadi_PEMag_14} \cite{steel_IEMag_10}. It is known for its high efficiency and power density.

High performance control of the PMSM can be achieved using vector control \cite{leonhard} \cite{chiasson} \cite{boldea_IEMag_08} \cite{betin_iemag_14}, which relies on the two-reactance theory developed by Park \cite{park_aiee_29} \cite{park_aiee_33}, where an accurate knowledge of the rotor position is required.

For many reasons, mainly for cost reduction and reliability increase \cite{pacas_iemag_11}, mechanical sensorless control of electrical drives has attracted the attention of researchers as well as many large manufacturers \cite{vas} \cite{holtz_tie_06} \cite{acarnley_tie_06} \cite{finch_tie_08}; mechanical sensors are to be replaced by algorithms that estimate the rotor speed and position, based on electrical sensors measurement.

One interesting sensorless technique is the observer-based one, which consists of sensing the machine currents and voltages, and using them as inputs to a state observer \cite{luenberger_tac_71} that estimates the rotor angular speed and position. There exists a tremendous variety of observers for PMSM in the literature \cite{benjak_icem_10_2}. Kalman filter \cite{kalman_asme_60} \cite{bolognani_tie_99} \cite{bolognani_tia_03} and sliding mode observers \cite{utkin_iecon_02} \cite{foo_tie_10} are among the most widely used observer algorithms in PMSM sensorless drives \cite{xu_tie_12}. Nevertheless, other nonlinear observers \cite{khalil} are also developed \cite{solsona_tie_96} \cite{zhu_tie_01} \cite{ortega_tcst_11} \cite{boussak_tpe_12}. Observer-based techniques rely on the machine mathematical model. Hence, depending on the modeling approach, three categories of these techniques can be distinguished:
\begin{itemize}
\item Electromechanical model-based observers \cite{dhaouadi_tpe_91} \cite{solsona_tie_96} \cite{bolognani_tie_99} \cite{zhu_tie_01} \cite{boussak_tpe_05} \cite{fadel_epe_07}.
\item Back electromotive force (EMF)-based observers \cite{chen_tie_03} \cite{mobarakeh_tia_04} \cite{akrad_tie_11} \cite{hijazi_sled_12}.
\item Flux-based observers \cite{ortega_tcst_11} \cite{koonlaboon_ias_05} \cite{boldea_tpe_08} \cite{foo_tpe_10}  \cite{koteich_ecmsm_13}.
\end{itemize}

Another sensorless technique is the high frequency injection (HFI) based technique \cite{corley_tia_98} \cite{arias_tie_06} \cite{diallo_tec_15}. Some authors propose to combine Observer-based and HFI techniques \cite{abry_epe_11} \cite{zgorski_sled_12} \cite{koteich_ecc_15}.

The main problem of the PMSM observer-based sensorless techniques is the deteriorated performance in low- and zero-speed operation conditions \cite{vaclavek_tie_13}. This problem is usually treated from observer's stability point of view, whereas the real problem remains hidden: it lies in the so-called \textit{observability conditions} of the machine. Over the past few years, a promising approach, based on the local weak observability concept \cite{hermann_tac_77}, has been used in order to better understand the deteriorated performance of sensorless drives.

Even though several papers have been published about the PMSM local observability, none of these papers presents well elaborated results for both salient and non-salient PMSMs, especially at standstill: 
\begin{itemize}
\item \textit{Surface-mounted PMSM (SPMSM)}: in \cite{zhu_tie_01} and \cite{ezzat_cdc_10} the SPMSM observability is studied; only the output first order derivatives are evaluated, and the conclusion is that the SPMSM  local observability cannot be guaranteed if the rotor speed is null. In \cite{zaltni_cdc_10} higher order derivatives of the output are investigated, and it is shown that the SPMSM can be locally observable at standstill if the rotor acceleration is nonzero.
\item \textit{Internal PMSM (IPMSM)}: concerning the IPMSM, the conclusions in \cite{zaltni_cdc_10} are unclear, i.e. no explicit \textit{practical} observability conditions are given. More interesting results are presented in \cite{vaclavek_tie_13}, where explicit conditions, expressed in the rotating reference frame, are presented. However, the analysis of the results in \cite{vaclavek_tie_13} remains unclear and yet inaccurate. In \cite{koteich_ecc_15}, a unified approach is adopted for synchronous machines observability study; the PMSM is treated as a special case of the generalized synchronous machine without further analysis.
\end{itemize}

The present paper is intended to investigate the PMSM observability for the electromechanical model-based observers. Both salient type (IPMSM) and non-salient type (SPMSM) machines observability analysis is detailed. A special attention is drawn to the standstill operation condition.

After this introduction, the paper is organized as follows: the local weak observability theory is presented in section II. Section III is dedicated to the PMSM's mathematical model in both stator and rotor reference frames. The observability analysis of the IPMSM is presented in section IV, whereas the SPMSM observability is analyzed in section V. Illustrative simulations are presented in section VI to validate the theoretical study. Conclusions are drawn in section VII.

\section{Local Observability Theory}

The \emph{local weak observability} concept \cite{hermann_tac_77}, based on the rank criterion, is presented in this section. The systems of the following form (denoted $\Sigma$) are considered:\\
\small
\begin{equation}
\Sigma :\left\{
\begin{aligned}
\dot{x} &= f\left(x(t), u(t)\right)\\
y &= h\left(x(t)\right)
\label{sigma}
\end{aligned}
\right.\end{equation}
\normalsize
where $x \in X \subset \mathbb{R}^n$ is the state vector, $u \in U \subset \mathbb{R}^m$ is the control vector (input), $y \in \mathbb{R}^p $ is the output vector, $f$ and $h$ are $C^\infty$ functions. The observation problem can be then formulated as follows \cite{besancon}:
\emph{Given a system described by a representation \eqref{sigma}, find an accurate estimate $\hat{x}(t)$ for $x(t)$ from the knowledge of $u(\tau)$, $y(\tau)$ for $0 \leq \tau \leq t$}.

\subsection{Observability rank condition}
The system $\Sigma$ is said to satisfy the observability rank condition at $x_0$ if the observability matrix, denoted by $\mathcal{O}_y(x)$, is full rank at $x_0$. $\mathcal{O}_y(x)$ is given by:
\small
\begin{equation}
\mathcal{O}_y(x) = \frac{\partial}{\partial x}\left[ \begin{matrix}
\mathcal{L}^0_fh(x) \\ 
\mathcal{L}_fh(x) \\ 
\mathcal{L}_f^2h(x) \\
\ldots \\ 
\mathcal{L}_f^{n-1}h(x)
\end{matrix} \right]_{x=x_0}^T
\end{equation}
\normalsize

where $\mathcal{L}_f^{k}h(x)$ is the $k$th-order \emph{Lie derivative} of the function $h$ with respect to the vector field $f$. 
It is given by:
\small
\begin{eqnarray}
 \mathcal{L}_f h(x) &=& \frac{\partial h(x)}{\partial x} f(x) \\ \mathcal{L}_f^k h(x) &=& \mathcal{L}_f \mathcal{L}_f^{k-1} h(x)\\
 \mathcal{L}_f^0 h(x) &=& h(x)
\end{eqnarray}
\normalsize

\subsection{Observability theorem}
A system $\Sigma$ \eqref{sigma} satisfying the observability rank condition at $x_0$ is locally weakly observable at $x_0$. More generally, a system $\Sigma$ \eqref{sigma} satisfying the observability rank condition for any $x_0$, is locally weakly observable. Rank criterion gives only a sufficient condition for local weak observability.

\section{PMSM Mathematical model}
Permanent magnet synchronous machines are electromechanical systems that can be mathematically represented using generalized Ohm's, Faraday's, and Newton's second Law. This section presents the PMSM model in two reference frames: the stator reference frame $\alpha \beta$, and the rotor reference frame $dq$ \cite{park_aiee_29} \cite{park_aiee_33}.

The assumption of linear lossless magnetic circuit is adopted, with sinusoidal distribution of the stator magnetomotive force (MMF). The machine parameters are considered to be known and constant. Nevertheless, the parameters variation does not call the observability study results into question; it impacts the observer performance, which is beyond the scope of this study.

\subsection{PMSM model in the stator reference frame}
The mathematical model of the PMSM in the stator reference frame can be written as:
\begin{equation}
\begin{aligned}
\frac{d\mathcal{I}}{dt}  &= {\mathfrak{L}^{-1}} \left( \mathcal{V}
-{\mathfrak{R}^{eq}}\mathcal{I}
-\psi_r \mathcal{C}'(\theta) \omega \right) \\
\frac{d\omega}{dt} &= \frac{p}{J} \left( T_m - T_l \right) \\
\frac{d\theta}{dt} &= \omega 
\end{aligned}\label{model_ipmsm}
\end{equation}
where $\mathcal{I} = \begin{bmatrix}
   {{i}_{\alpha }}  &
   {{i}_{\beta }}
\end{bmatrix}^T$ and $ \mathcal{V} = \begin{bmatrix}
   {{v}_{\alpha }} & {{v}_{\beta }}  
\end{bmatrix}^T$ stand for currents and voltages vectors in the $\alpha \beta$ reference frame, $\omega$ is the electrical speed of the rotor, $\theta$ is its electrical position, $\psi_r$ is the rotor permanent magnet flux. $\mathfrak{L}$ is the inductance matrix, $\mathfrak{R}_{eq}$ is the equivalent resistance matrix:
\begin{eqnarray}
\mathfrak{R}^{eq} = \mathfrak{R} + \omega \mathfrak{L}' =  \begin{bmatrix}
R & 0\\ 0 & R
\end{bmatrix} + \omega  \frac{\partial \mathfrak{L}}{\partial \theta} 
\end{eqnarray}
$R$ is the resistance of one stator winding. $p$ is the number of pole pairs, $J$ is the inertia of the rotor with the load, $T_l$ is the resistant torque, and $T_m$ is the motor torque. $\mathcal{C}'(\theta)$ denotes the partial derivative of $\mathcal{C}(\theta)$ with respect to $\theta$:
\begin{eqnarray}
\mathcal{C}(\theta) = \begin{bmatrix}
\cos \theta \\ \sin \theta
\end{bmatrix} ~~;~~ \mathcal{C}'(\theta) = \frac{\partial \mathcal{C}(\theta)}{\partial \theta} = \begin{bmatrix}
-\sin \theta \\ \cos \theta
\end{bmatrix}
\end{eqnarray}

The model \eqref{model_ipmsm} can be fitted to the structure $\Sigma$ \eqref{sigma} by taking:
\begin{eqnarray}
x = \begin{bmatrix}
\mathcal{I}^T & \omega & \theta
\end{bmatrix}^T
~~;~~
y = \mathcal{I}~~;~~
u = \mathcal{V}\\
f(x,u) = \begin{bmatrix}
\frac{d\mathcal{I}^T}{dt} & \frac{d\omega}{dt} & \frac{d\theta}{dt}
\end{bmatrix}^T
~~;~~
h(x) = \mathcal{I}
\label{ss_vectors}
\end{eqnarray}

\subsubsection{IPMSM}
\begin{figure}[!t]
\centering
\includegraphics[scale=0.75]{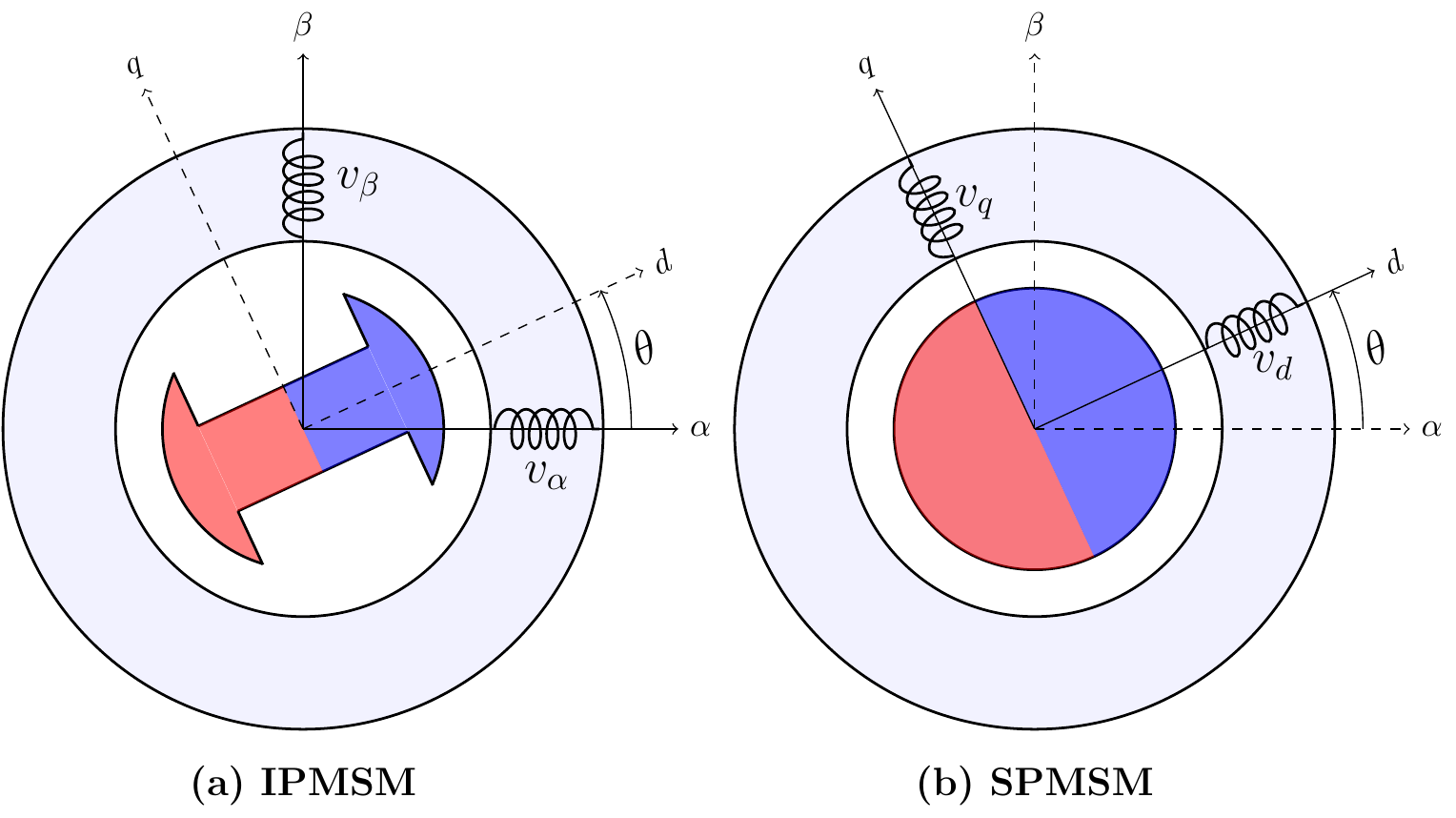}
\caption{Schematic representation of IPMSM and SPMSM}
\label{ipmsm}
\end{figure}
The IPMSM is a salient rotor machine, then its inductance matrix $\mathfrak{L}$ is a position-dependent matrix:
\begin{eqnarray}
\mathfrak{L}=\left[ \begin{matrix}
   {{L}_{0}}+{{L}_{2}}\cos 2\theta  & {{L}_{2}}\sin 2\theta \\
   {{L}_{2}}\sin 2\theta  & {{L}_{0}}-{{L}_{2}}\cos 2\theta
\end{matrix} \right] 
\end{eqnarray}
where $L_0$ and $L_2$ are the average and differential spatial inductances. The IPMSM produced torque is:
\begin{eqnarray}
T_m &=& \frac{3p}{2} \left[\psi_r (i_{\beta} \cos \theta - i_{\alpha} \sin \theta)\right. \\
&& \left. - L_2 \left((i_{\alpha}^2 - i_{\beta}^2 ) \sin 2\theta - 2 i_{\alpha} i_{\beta} \cos 2\theta \right) \right] \nonumber
\end{eqnarray}

\subsubsection{SPMSM}
The SPMSM is a non-salient rotor machine, its model can be derived from the IPMSM one by assuming $L_2$ to be null:
\begin{eqnarray}
L_2 = 0
\end{eqnarray}

\subsection{PMSM model in the rotor reference frame}
Stator currents and voltages in the $dq$ rotating reference frame (Fig. \ref{ipmsm}) are calculated from those in the $\alpha \beta$ reference frame using the Park transform:
\begin{eqnarray}
\mathcal{X}_{dq} = \mathcal{P}{^{-1}}(\theta) \mathcal{X}_{\alpha \beta}
\label{park}
\end{eqnarray}
where
\begin{eqnarray}
\mathcal{P}(\theta) = \begin{bmatrix}
\cos \theta & - \sin \theta \\
\sin \theta & \cos \theta
\end{bmatrix}
\end{eqnarray}
The vector $\mathcal{X}$ stands for currents or voltages vector. The mathematical model of the PMSM in the rotor reference frame can be then written under the general form \eqref{model_ipmsm}:
\begin{equation}
\begin{aligned}
\frac{d\mathcal{I}_{dq}}{dt}  &= {\mathfrak{L}_{dq}^{-1}} \left( \mathcal{V}_{dq}
-{\mathfrak{R}_{dq}^{eq}}\mathcal{I}_{dq}
-\psi_r \mathcal{C}'(0) \omega \right) \\
\frac{d\omega}{dt} &= \frac{p}{J} \left( T_m - T_l \right) \\
\frac{d\theta}{dt} &= \omega 
\end{aligned}\label{model_ipmsm_dq}
\end{equation}
where $\mathcal{I}_{dq} = \begin{bmatrix}
   {{i}_{d }}  &
   {{i}_{q }}
\end{bmatrix}^T$ and $ \mathcal{V}_{dq} = \begin{bmatrix}
   {{v}_{d }} & {{v}_{q }}  
\end{bmatrix}^T$ stand for currents and voltages vectors in the $dq$ reference frame. Inductance and equivalent resistance matrices can be written as:
\begin{eqnarray}
\mathfrak{L}_{dq} = \begin{bmatrix}
L_d & 0\\ 0 & L_q
\end{bmatrix} ~~;~~ \mathfrak{R}_{dq}^{eq} = \mathfrak{R} + \omega \mathbb{J}_2 \mathfrak{L}_{dq}
\end{eqnarray}
where
\begin{eqnarray}
\mathbb{J}_2 = \mathcal{P}\left(\frac{\pi}{2}\right) = \begin{bmatrix}
0 & -1 \\1 & 0
\end{bmatrix}
\end{eqnarray}
$L_d$ denotes the (direct) $d-$axis inductance, and $L_q$ denotes the (quadrature) $q-$axis inductance:
\begin{eqnarray}
L_d &=& L_0 + L_2\\
L_q &=& L_0 - L_2
\end{eqnarray}

The motor torque can be written as:
\begin{eqnarray}
T_m = \frac{3}{2} p \left(L_{\delta} i_d + \psi_r \right) i_q
\end{eqnarray}
with
\begin{eqnarray}
L_\delta = L_d - L_q = 2 L_2
\end{eqnarray}

The above $dq$ model is valid for the IPMSM ($L_d \neq L_q$). The SPMSM model can be derived using the following equations:
\begin{eqnarray}
L_d = L_q = L_0 \implies L_\delta = 0
\end{eqnarray}

\section{IPMSM Observability}
Observability of the system \eqref{model_ipmsm} is studied in the sequel. The system \eqref{model_ipmsm} is a $4^{th}$ order system. Its observability matrix should contain the gradient of the output and its derivatives up to the $3^{rd}$ order. In this section, only the first order derivatives are calculated, higher order derivatives are very difficult to evaluate and to deal with. The ``partial'' observability matrix is:
\begin{eqnarray}
\mathcal{O}_{y1} 
= \frac{\partial (y, \dot{y})}{\partial x} = \begin{bmatrix}
\mathbb{I}_{2} & \mathbb{O}_{2 \times 1} & \mathbb{O}_{2 \times 1} \\
\frac{\partial}{\partial \mathcal{I} }\left(\frac{d\mathcal{I}}{dt} \right) & \frac{\partial}{\partial \omega }\left(\frac{d\mathcal{I}}{dt} \right) & \frac{\partial}{\partial {\theta} }\left(\frac{d\mathcal{I}}{dt} \right)
\label{obsv_matrix_ipmsm}
\end{bmatrix} 
\end{eqnarray}
where $\mathbb{I}_{n}$ is the $n \times n$  identity matrix, and $\mathbb{O}_{n \times m}$ is an $n \times m$ zero matrix, and:
\begin{eqnarray}
\frac{\partial}{\partial \mathcal{I} }\left(\frac{d\mathcal{I}}{dt} \right) &=& -\mathfrak{L}^{-1} \mathfrak{R}^{eq} \nonumber\\
\frac{\partial}{\partial \omega }\left(\frac{d\mathcal{I}}{dt} \right) &=& -\mathfrak{L}^{-1} \left(\mathfrak{L}' \mathcal{I} + \psi_r \mathcal{C}'(\theta) \right)\\
\frac{\partial}{\partial {\theta} }\left(\frac{d\mathcal{I}}{dt} \right) &=& (\mathfrak{L}^{-1})' \mathfrak{L} \frac{d \mathcal{I}}{dt} - \mathfrak{L}^{-1} \left(\mathfrak{L}'' \mathcal{I} - \psi_r \mathcal{C}(\theta) \right) \omega \nonumber
\end{eqnarray}
$\mathfrak{L}'$ and $\mathfrak{L}''$ denote, respectively, the first and second partial derivatives of $\mathfrak{L}$ with respect to $\theta$:
\begin{eqnarray}
\mathfrak{L}' = \frac{\partial}{\partial \theta} \mathfrak{L}~~~;~~
\mathfrak{L}'' = \frac{\partial}{\partial \theta} \mathfrak{L}'
\end{eqnarray}
The determinant $\Delta_{y1} $ of the sub-matrix  \eqref{obsv_matrix_ipmsm} is calculated using symbolic math software. In order to make the interpretation easier, the determinant is expressed in the rotating $dq$ reference frame using the equation \eqref{park}.

The determinant $\Delta_{y1} $ is given by:
\begin{eqnarray}
\Delta_{y1} &=&\frac{1}{L_d L_q} \left[\left(L_\delta i_d + \psi_r \right)^2 + L_\delta^2 i_q^2
\right]\omega \nonumber \\
&&+~ \frac{L_\delta}{L_d L_q} \left[
 L_\delta\frac{di_d}{dt} i_q - 
\left(L_\delta i_d + \psi_r \right)  \frac{di_q}{dt}
\right]
\end{eqnarray}

The observability condition $\Delta_{y1} \neq 0$ can be written as:
\begin{eqnarray}
\omega \neq \frac{\left(L_\delta i_d + \psi_r \right) L_\delta \frac{di_q}{dt} - L_\delta \frac{di_d}{dt} L_\delta i_q}{\left(L_\delta i_d + \psi_r \right)^2 + L_\delta^2 i_q^2}
\end{eqnarray}
which gives:
\begin{eqnarray}
\omega \neq   \frac{d}{dt} \arctan \left(\frac{L_\delta i_q}{L_\delta i_d + \psi_r}\right)
\label{cond_w}
\end{eqnarray}

The equation \eqref{cond_w} defines a fictitious \textit{observability vector}, denoted by $\Psi_\mathcal{O}$, that has the following components in the $dq$ reference frame:
\begin{eqnarray}
\Psi_{\mathcal{O}d} &=& L_\delta i_d + \psi_r\\
\Psi_{\mathcal{O}q} &=& L_\delta i_q
\end{eqnarray}

\begin{figure}[!t]
\centering
\includegraphics[scale=1.2]{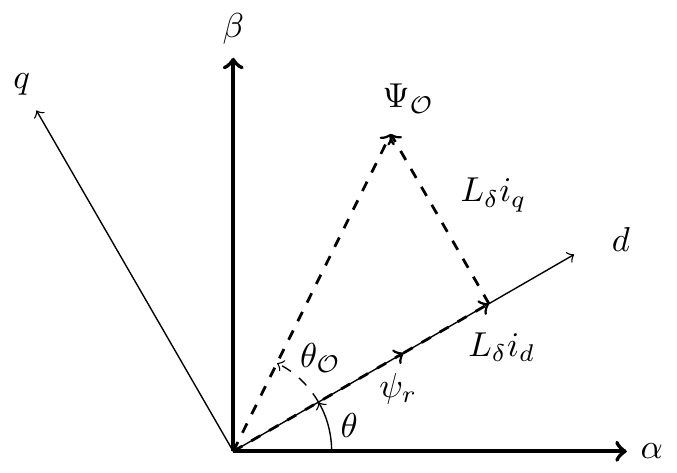}
\caption{Vector diagram of the fictitious observability vector (dashed)}
\label{obsv_vector}
\end{figure}

Then, the condition \eqref{cond_w} can be formulated as:
\begin{eqnarray}
\omega \neq   \frac{d}{dt} \theta_\mathcal{O}
\label{condition}
\end{eqnarray}
where $\theta_\mathcal{O}$ is the phase of the vector $\Psi_\mathcal{O}$ in the rotating reference frame (Fig. \ref{obsv_vector}). The following sufficient condition for the PMSM local observability can be stated: the rotational speed of the fictitious vector $\Psi_\mathcal{O}$ in the $dq$ reference frame should be different from the rotor electrical angular speed in the stator reference frame. At standstill, the above condition becomes: the vector $\Psi_\mathcal{O}$ should keep changing its orientation in order to ensure the local observability. 

It turns out that the $d-$axis component of the vector $\Psi_\mathcal{O}$ is nothing but the so-called ``active flux'' introduced by Boldea \textit{et al.} in \cite{boldea_tpe_08} (also called ``fictitious flux'' by Koonlaboon \textit{et al.} \cite{koonlaboon_ias_05}), which is, by definition, the torque producing flux aligned to the rotor $d-$axis. The $q-$axis component of the vector $\Psi_\mathcal{O}$ is related to the saliency ($L_\delta$) of the machine.

Some authors \cite{vaclavek_tie_13} present conclusions that stator current space vector should change not only its magnitude, but also direction, in the rotating reference frame in order to ensure motor observability at standstill. However, Fig.~\ref{obsv_vector} shows that the stator current space vector can change both its magnitude and direction without fulfilling the condition \eqref{condition}. 

For the SPMSM ($L_\delta = 0$), the fictitious observability vector is equal to the rotor PM flux vector, which is fixed in the $dq$ reference frame. This means that the observability problem arises only at standstill, as shown in the next section.

\section{SPMSM observability}
Fortunately, the SPMSM model is less complex than the IPMSM one, which makes the investigation of higher output derivatives possible.
In this section, the electromechanical model observability is studied. Furthermore, thanks to the simplicity of the SPMSM equations, the observabiltiy of two other models, namely the back-emf and flux models, is studied. 

\subsection{Electromechanical model observability}
In the case of SPMSM, the observability matrix $\mathcal{O}_{y}$ can be evaluated up to the $3^{rd}$ order output derivatives. $\mathcal{O}_y$ is an $8 \times 4$ matrix. There are 70 possible $4 \times 4$ sub-matrices. For convenience, the first two lines are always taken, together with lines that correspond to the same derivation order. This reduces the choices to the following 3 possible sub-matrices: 
\begin{itemize}
\item $\mathcal{O}_{y1}$, which includes the first 4 lines of $\mathcal{O}_y$:
\begin{eqnarray}
\mathcal{O}_{y1} = \begin{bmatrix}
1 & 0 & 0 & 0 \\
0 & 1 & 0 & 0 \\
-\frac{R}{L_0} & 0 & \frac{\psi_r}{L_0} \sin \theta & \omega \frac{\psi_r}{L_0} \cos \theta\\
0 & -\frac{R}{L_0} & -\frac{\psi_r}{L_0} \cos \theta & \omega \frac{\psi_r}{L_0} \sin \theta
\end{bmatrix}
\end{eqnarray}
its determinant is
\begin{eqnarray}
\Delta_{y1} = \omega\left(\frac{\psi_r}{L_0}\right)^2
\label{det_spmsm_1}
\end{eqnarray}
Thus, the local observability is guaranteed if the rotor speed is nonzero, but not in the case of zero speed.

\item $\mathcal{O}_{y2}$, which includes the $1^{st}$, $2^{nd}$, $5^{th}$ and $6^{th}$ lines of $\mathcal{O}_y$, its determinant is :
\begin{eqnarray}
\Delta_{y2} = \frac{\psi_r^2}{L_0^2} \left[ \left(2\omega^2 +  \frac{R^2}{L_0^2}   + \frac{3 p^2}{J} \psi_r i_d\right) \omega  -  \frac{R}{L_0}  \frac{d\omega}{dt} \right]
\end{eqnarray}
Thus, in the case of zero speed operation, the SPMSM is observable if the acceleration is different than zero ($\dot{\omega} \neq 0$);  this corresponds to the case where the motor changes its rotation direction.

\item $\mathcal{O}_{y3}$, which includes the $1^{st}$, $2^{nd}$, $7^{th}$ and $8^{th}$ lines of $\mathcal{O}_y$. Its determinant $\Delta_{y3}$ cannot be written because it is lengthy. Nevertheless, substituting the rank deficiency conditions of the sub-matrix $\mathcal{O}_{y2}$ ($\omega = 0$ and $\dot{\omega}=0$) in $\Delta_{y3}$, under the assumption of very slow resistant torque variation ($\dot{T_l}=0$), gives:
\begin{eqnarray}
\Delta_{y3}|_{\Delta_{y2}=0} =  \frac{\psi_r^2}{L_0^2} \left[ \frac{R^2}{L_0^2} - \frac{3p^2}{2J}  (L_0 i_d + \psi_r) \frac{\psi_r}{L_0} \right] \frac{d^2\omega}{dt^2} 
\end{eqnarray}
where
\begin{eqnarray}
\frac{d^2\omega}{dt^2} = \frac{3p^2}{2J} \psi_r \frac{di_q}{dt}
\end{eqnarray}
\end{itemize}

If the speed is identically zero ($\omega \equiv 0$), the SPMSM model reduces to:
\begin{eqnarray}
\frac{d \mathcal{I}}{dt} &=& \frac{1}{L_0} \left(\mathcal{V} - R \mathcal{I}\right) \nonumber\\
\frac{d\omega}{dt} &=& 0 \label{model_spmsm_w0}\\
\frac{d \theta}{dt} &=& 0 \nonumber
\end{eqnarray}
and the output derivatives are:
\begin{eqnarray}
\frac{d \mathcal{I}}{dt} &=& \frac{1}{L_0} \left(\mathcal{V} - R \mathcal{I}\right) \label{deriv_1}\\
\frac{d^2 \mathcal{I}}{dt^2} &=& \frac{1}{L_0} \left(\frac{d \mathcal{V}}{dt} - R \frac{d \mathcal{I}}{dt}\right)\\
\vdots \nonumber\\
\frac{d^{n+1} \mathcal{I}}{dt^{n+1}} &=& \frac{1}{L_0} \left(\frac{d^{n} \mathcal{V}}{dt^{n}} - R \frac{d^{n} \mathcal{I}}{dt^{n}}\right) \label{deriv_n}
\end{eqnarray}
In this case, the observability matrix is:
\begin{eqnarray}
\mathcal{O}_y|_{\omega \equiv 0} = \begin{bmatrix}
1 & 0  & 0 & 0\\
0 & 1  & 0 & 0\\
-\frac{R}{L_0} & 0 & \frac{\psi_r}{L_0} \sin \theta & 0\\
0 & -\frac{R}{L_0} & -\frac{\psi_r}{L_0} \cos \theta & 0\\
(-\frac{R}{L_0})^2 & 0  & -\frac{R}{L_0}\frac{\psi_r}{L_0} \sin \theta & 0\\
0 & (-\frac{R}{L_0})^2  & \frac{R}{L_0}\frac{\psi_r}{L_0} \cos \theta & 0 \\
(-\frac{R}{L_0})^3 & 0  & (-\frac{R}{L_0})^2\frac{\psi_r}{L_0} \sin \theta & 0\\
0 & (-\frac{R}{L_0})^3  & -(-\frac{R}{L_0})^2\frac{\psi_r}{L_0} \cos \theta  & 0
\end{bmatrix}
\label{obsv_matrix_w0}
\end{eqnarray}
The following recurrence can be obtained from \eqref{deriv_1}-\eqref{obsv_matrix_w0} for higher dimension observability matrices:
\begin{eqnarray}
\frac{\partial}{\partial x} \mathcal{L}_f^k h = -\frac{R}{L_0} \frac{\partial}{\partial x} \mathcal{L}_f^{k-1}h |_{\omega \equiv 0}
\end{eqnarray}
Therefore, even if higher order derivatives are evaluated, no additional information about the rotor position can be extracted. Hence, the standstill operation condition presents a singularity from observability viewpoint. 

Physically speaking, if the non-salient (cylindrical) rotor is fixed with respect to the stator windings, it will have no effect on the electromagnetic behaviour of these windings, and its position cannot be identified with the model \eqref{model_spmsm_w0}. 

\begin{figure}[!t]
\centering
\includegraphics[scale=1.2]{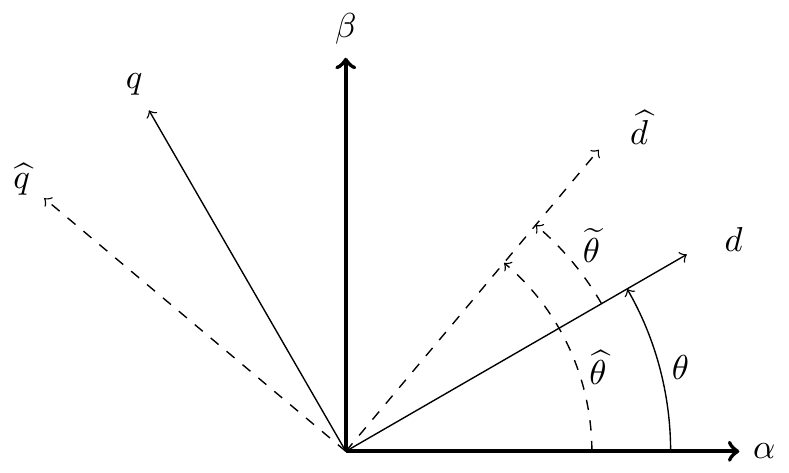}
\caption{Vector diagram of stator (thick), estimated (dashed) and real rotor reference frames}
\label{inj_error}
\end{figure}

One solution of this problem is proposed in \cite{abry_epe_11}. It combines HFI technique with a state observer algorithm; a sinusoidal voltage is injected on the direct ($\widehat{d}-$) axis in the estimated rotating reference frame (Fig. \ref{inj_error}), which results in a vibration of the rotor (nonzero speed) only if the position is not correctly estimated. It is proved in \cite{abry_epe_11} that if the following HF voltage 
\begin{eqnarray} 
v_{\hat{d}} = V_{hf} \cos (\omega_{hf} t)
\end{eqnarray}
is injected in the estimated $\hat{d}\hat{q}$ rotating reference frame (obtained by the Park transformation using the estimated position $\hat{\theta}$), the determinant $\Delta_{y1}$ \eqref{det_spmsm_1} becomes:
\begin{eqnarray}
\Delta_{y1} = - \frac{\psi_r^2}{L_0^2} \omega + \frac{\psi_r}{L_0^2} V_{hf} \cos (\omega_{hf} t) \sin \widetilde{\theta }
\end{eqnarray} 
where $\widetilde{\theta}$ stands for the position estimation error. If the position is not correctly estimated ($\widetilde{\theta} \neq 0$), the local observability is guaranteed at standstill.

Another solution is proposed in \cite{scaglione_icems_11} \cite{scaglione_tie_12}. It consists of adding a position-dependent source, $g(\theta)$, to the available measurements. The origin of $g(\theta)$ is found in the stator iron local B-H hysteresis loops. It is shown that this signal is highly position-dependent, and it can be approximated by a linear function:
\begin{eqnarray}
g(\theta) = a \theta + b
\end{eqnarray}
where $a$ and $b$ are two constants. The new output vector becomes: \begin{eqnarray}
y = \begin{bmatrix}
i_\alpha & i_\beta & a\theta + b
\end{bmatrix}^T
\end{eqnarray}
Therefore, even at standstill, the position is observable.

\subsection{Back-EMF model observability}
The state-space model for the back-EMF based observer can be written as:
\begin{eqnarray}
\frac{d \mathcal{I}}{dt} &=& \frac{1}{L_0} \left( \mathcal{V} - \mathfrak{R} \mathcal{I} - \mathcal{E}\right)\\
\frac{d \mathcal{E}}{dt} &=& \left(\frac{\dot{\omega}}{\omega} \mathbb{I}_2 + \omega \mathbb{J}_2 \right) \mathcal{E}
\end{eqnarray}
where $\mathcal{E} = \left[ e_\alpha ~~ e_\beta\right]^T$ stands for the back-EMF vector in the stator $\alpha \beta$ reference frame:
\begin{eqnarray}
e_\alpha &=& - \omega \psi_r \sin \theta \label{e_a} \\
e_\beta &=& \omega \psi_r \cos \theta \label{e_b}
\end{eqnarray}
Then the rotor speed and position can be calculated from the back-EMF components using the following relationships:
\begin{eqnarray}
\theta &=& \arctan\left(-\frac{e_\alpha}{e_\beta}\right) \label{theta_e_spmsm}\\
\omega &=& \frac{1}{\psi_r} \sqrt{e_\alpha^2 + e_\beta^2} \label{omega_e_spmsm}
\end{eqnarray}

The observability analysis is done for the following state, input and output vectors:
\begin{eqnarray}
x = \begin{bmatrix}
\mathcal{I}^T & \mathcal{E}^T
\end{bmatrix}^T ~~;~~
u = \mathcal{V} ~~;~~
y = \mathcal{I}
\end{eqnarray}

The sub-matrix made of the first 4 lines of the observability matrix is studied:
\begin{eqnarray}
\mathcal{O}_{y_1}  = \begin{bmatrix}
1 & 0 & 0 & 0\\
0 & 1 & 0 & 0\\
-\frac{R}{L_0} & 0 & -\frac{1}{L_0} & 0\\
0 & -\frac{R}{L_0} & 0 & -\frac{1}{L_0}
\end{bmatrix}
\end{eqnarray}
its determinant is.
\begin{eqnarray}
\Delta_{y1} = \frac{1}{L_0^2}
\end{eqnarray}

This implies that the system is observable even at zero speed (no need to calculate higher order derivatives). However, at standstill, the back-EMF components are both null (see equations \eqref{e_a} \eqref{e_b}), and the position \eqref{theta_e_spmsm} is indeterminate. The problem remains the same: the rotor position is not observable at standstill.

\subsection{Flux-based model observability}
The state-space model for the flux-based observer can be written as:
\begin{eqnarray}
\frac{d \mathcal{I}}{dt} &=& \frac{1}{L_0} \left( \mathcal{V} - \mathfrak{R} \mathcal{I} - \omega \mathbb{J}_2 \Psi_r\right)\\
\frac{d \Psi_r}{dt} &=& \omega \mathbb{J}_2 \Psi_r
\end{eqnarray}
where $\Psi_r = \left[ \psi_{r_\alpha} ~~ \psi_{r_\beta} \right]^T$ stands for the rotor magnetic flux vector in the stator $\alpha \beta$ reference frame:
\begin{eqnarray}
\psi_{r_\alpha} &=& \psi_r \cos \theta \\
\psi_{r_\beta} &=& \psi_r \sin \theta
\end{eqnarray}
The rotor position is given by:
\begin{eqnarray}
\theta = \arctan\left(\frac{\psi_{r_\beta}}{\psi_{r_\alpha}}\right)
\end{eqnarray}

The observability analysis is done for the following state, input and output vectors:
\begin{eqnarray}
x = \begin{bmatrix}
\mathcal{I}^T & \Psi_r^T
\end{bmatrix}^T ~~;~~
u = \mathcal{V} ~~;~~
y = \mathcal{I}
\end{eqnarray}

The following determinants correspond respectively to the $1^{st}$, $2^{nd}$ and $3^{rd}$ order output derivatives:
\begin{eqnarray}
\Delta_{y1} &=& \frac{\omega^2}{L_0^2}\\
\Delta_{y2} &=& \frac{\omega^2}{L_0^4} \left(R^2 + L_0^2\omega^2\right)\\
\Delta_{y3} &=& \frac{\omega^2}{L_0^6} \left(R^4 + L_0^4 \omega^4 - R^2 L_0^2 \omega^2\right)
\end{eqnarray}
It is obvious that the system is not observable at standstill. The authors in \cite{zgorski_sled_12} propose an HF injection-based solution combined with a state observer.

\section{Illustrative simulations}
\begin{figure}[!t]
\begin{center}
\includegraphics[scale=1]{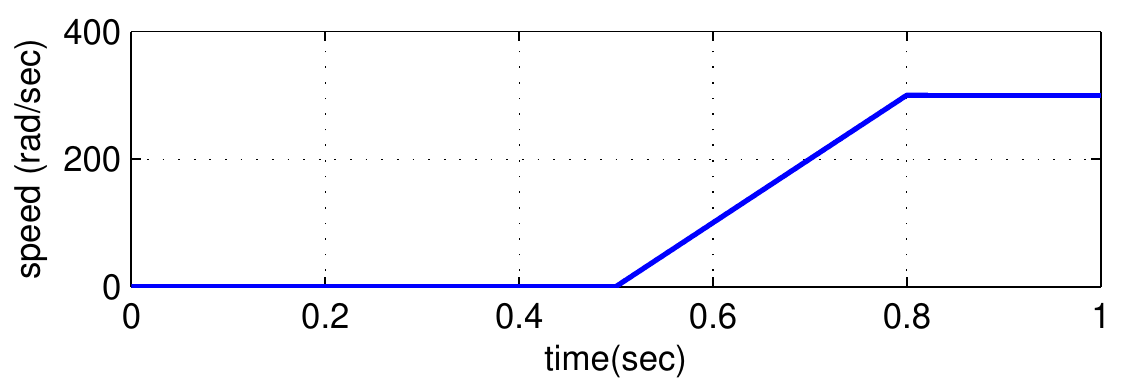}
\end{center}
\caption{Rotor speed profile}
\label{speed_profile}
\end{figure}

The present section is aimed at illustrating the previous observability analysis using numerical simulations. For this purpose, an extended Kalman filter (EKF) is designed. In order to make the study of some critical situations easier, the following operation mode is installed: the rotor position is considered to be driven by an external mechanical system, which imposes the speed profile shown in Fig.~\ref{speed_profile}. The currents are regulated, using standard proportional-integral (PI) controllers, to fit with the following set-points:
\begin{eqnarray}
\begin{array}{c c c c c}
i_d^* = 0~A & ; & i_q^* = 15~A
\end{array}
\end{eqnarray}

Both IPMSM and SPMSM are studied in the same simulation environment. The same machine parameters are used for both machines; the only difference is in the inductance $L_2$, which is null in the case of SPMSM (no saliency). The following HF current is added to the current $i_q$ during the time interval $[0.2~s., 0.5~s.]$:
\begin{eqnarray}
i_{q_{HF}} = 0.5 \sin 1000 \pi t ~~ A
\end{eqnarray}

\begin{table}[!t]
\renewcommand{\arraystretch}{1.3}
\label{table}
\caption{IPMSM Parameters}
\centering
\begin{tabular}{|l|c|}
\hline
Parameters & Value [Unit]\\
\hline
\hline
Number of pole pairs (${p}$) & 2\\
Stator resistance $R_s$ & $0.01$ [$\Omega$]\\
Direct inductance $L_d$ & $0.5$ [$mH$]\\
Quadratic inductance $L_q$ & $0.8$ [$mH$]\\
Rotor magnetic flux $\psi_r$ & $0.0225$ [$V.s/rad$]\\
\hline
\end{tabular}
\label{param}
\end{table}

The purpose is to compare the observer behaviour for both machines at standstill, with and without signal injection. The observer is operating in open-loop, the real position is fed to the controller in order to avoid stability issues in the observability analysis. Table \ref{param} shows the machine parameters.

\subsection{Extended Kalman Filter}
The EKF algorithm is described below:
\paragraph{Model linearization}
\begin{eqnarray}
A_k = \left.\frac{\partial f(x,u)}{\partial x}\right|_{x_k, u_k};~~~
C_k = \left.\frac{\partial h(x)}{\partial x}\right|_{x_k} 
\end{eqnarray}
\paragraph{Prediction}
\begin{eqnarray}
\hat{x}_{k+1/k} &=& \hat{x}_{k/k} + T_s f(\hat{x}_{k/k},u_k)\\
P_{k+1/k} &=& P_k + T_s(A_k P_k + P_k A_k^T) + Q_k
\end{eqnarray}
\paragraph{Gain}
\begin{equation}
K_k = P_{k+1/k} C_k^T(C_k P_{k+1/k} C_k^T + R_k)^{-1}
\end{equation}
\paragraph{Innovation}
\begin{eqnarray}
\hat{x}_{k+1/k+1} &=& \hat{x}_{k+1/k} + K_k(y - h(\hat{x}_{k+1/k}))\\
P_{k+1/k+1} &=& P_{k+1/k} - K_k C_k P_{k+1/k}
\end{eqnarray}
where $T_s$ is the sampling period.
\paragraph{Tuning} EKF tuning is done by the choice of covariance matrices $Q_k$ and $R_k$. In this work the following matrices are used:
\begin{eqnarray}
Q_k = \left[
\begin{matrix}
1 & 0 & 0 & 0\\
0 & 1 & 0 & 0\\
0 & 0 & 10^3 & 0\\
0 & 0 & 0 & 0.1
\end{matrix}
\right]~;~
R_k = \begin{bmatrix}
1 & 0\\
0 & 1
\end{bmatrix}
\end{eqnarray}
The Kalman filter tuning has an impact on the estimation dynamics, which is beyond the scope of this paper. The same EKF is applied to both salient and non-salient machines, in order to compare the position estimation at standstill under the same conditions.

\subsection{Position and Speed estimation}
\begin{figure}[!t]
\begin{center}
\includegraphics[scale=1]{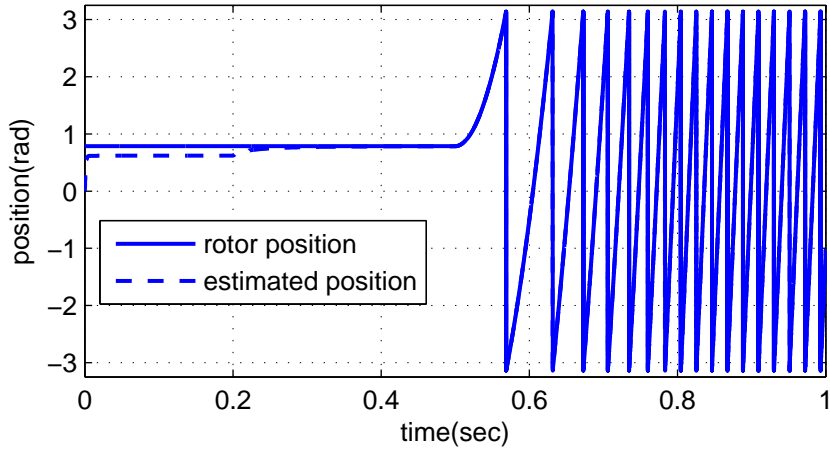}
\end{center}
\caption{Rotor real and estimated position of the IPMSM}
\label{position_ipmsm}
\end{figure}

\begin{figure}[!t]
\begin{center}
\includegraphics[scale=1]{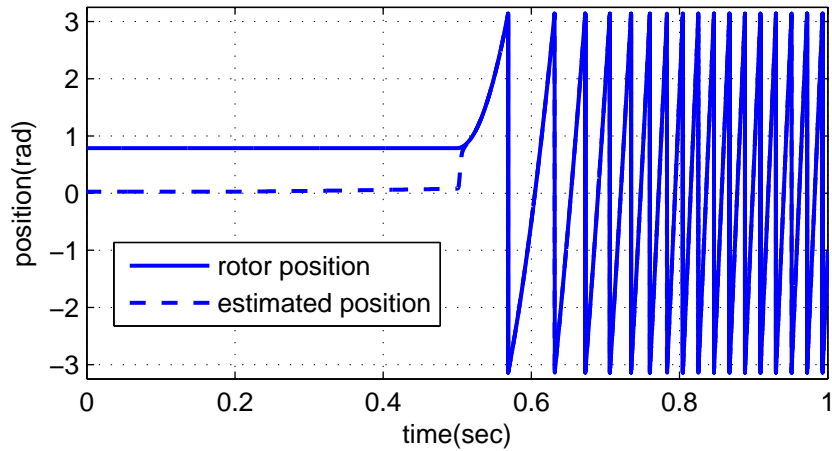}
\end{center}
\caption{Rotor real and estimated position of the SPMSM}
\label{position_spmsm}
\end{figure}

The EKF is initialized with a position error of $-\pi/4$. Fig.~\ref{position_ipmsm} and Fig.~\ref{position_spmsm} show respectively the position estimation for the IPMSM and SPMSM. At standstill, the following observations can be made:
\begin{itemize}
\item Before injecting the HF current, the position estimation of the IPMSM is more accurate than the SPMSM one.
\item After injecting the HF current, the IPMSM estimated position converges to the real position, whereas the SPMSM one slightly varies.
\item The SPMSM estimated position converges to the real position value as soon as the rotor accelerates.
\end{itemize}
These results are consistent with the observability study results; the IPMSM can be observable at standstill, whereas observability of the SPMSM cannot be guaranteed unless the rotor moves.


\begin{figure}[!t]
\begin{center}
\includegraphics[scale=1]{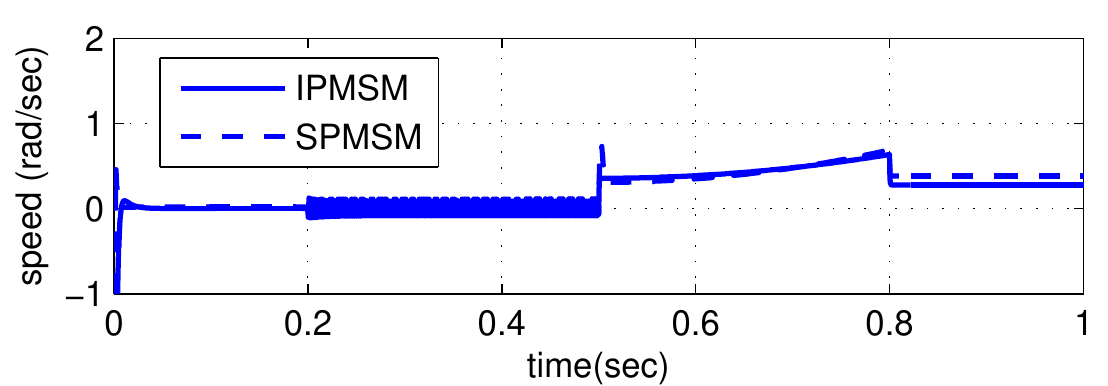}
\end{center}
\caption{Rotor speed estimation error}
\label{speed_error}
\end{figure}

The speed estimation error is shown in Fig.~\ref{speed_error}, it is almost the same for both machines.

\section{Conclusions}
The local observability study of the PMSM resulted in the definition of a fictitious \emph{observability vector}; the rotational speed of this observability vector in the rotor reference frame should be different from the rotor electrical speed in the stator reference frame to ensure the machine observability.

The results presented in this paper are valid for a wide range of brushless synchronous machines under the assumption of sinusoidal stator MMF distribution: PM synchronous, Brushless DC, PM stepper and PM assisted reluctance machines. Furthermore, if the rotor PM flux is considered to be zero, the results can be extended to synchonous reluctance machines.

\bibliographystyle{plain}
\bibliography{bib_SLED}

\end{document}